\documentclass{tran-l}
\usepackage{amsfonts}
\usepackage{amssymb}
\usepackage{xypic}

\begin{document}

\newtheorem{theorem}{Theorem}[subsection]
\newtheorem{lem}[theorem]{Lemma}
\newtheorem{cor}[theorem]{Corollary}
\newtheorem{prop}[theorem]{Proposition}

\theoremstyle{definition}
\newtheorem{definition}[theorem]{Definition}
\newtheorem{example}[theorem]{Example}
\newtheorem{xca}[theorem]{Exercise}

\theoremstyle{remark}
\newtheorem{remark}[theorem]{Remark}

\theoremstyle{plain}

\newtheorem{stheorem}{Theorem}[section]
\newtheorem{slem}[stheorem]{Lemma}
\newtheorem{scor}[stheorem]{Corollary}
\newtheorem{sprop}[stheorem]{Proposition}

\theoremstyle{definition}
\newtheorem{sdefinition}[stheorem]{Definition}
\newtheorem{sexample}[stheorem]{Example}
\newtheorem{sxca}[stheorem]{Exercise}

\theoremstyle{remark}
\newtheorem{sremark}[stheorem]{Remark}

\theoremstyle{plain}

%\numberwithin{equation}{section}
\def\limind{\mathop{\oalign{lim\cr
\hidewidth$\longrightarrow$\hidewidth\cr}}}
\def\boxit#1#2{\setbox1=\hbox{\kern#1{#2}\kern#1}%
\dimen1=\ht1 \advance\dimen1 by #1
\dimen2=\dp1 \advance\dimen2 by #1
\setbox1=\hbox{\vrule height\dimen1 depth\dimen2\box1\vrule}%
\setbox1=\vbox{\hrule\box1\hrule}%
\advance\dimen1 by .4pt \ht1=\dimen1
\advance\dimen2 by .4pt \dp1=\dimen2 \box1\relax}
\def \pext{\, \hbox{\boxit{0pt}{$\times$}}\,}

\newcommand{\sur}[2]{\genfrac{}{}{0pt}{}{#1}{#2}}

\let\cal\mathcal
\let\got\mathfrak
\def\AA{{\mathbf A}}
\def\BB{{\mathbf B}}
\def\CC{{\mathbf C}}
\def\DD{{\mathbf D}}
\def\EE{{\mathbf E}}
\def\FF{{\mathbf F}}
\def\GG{{\mathbf G}}
\def\HH{{\mathbf H}}
\def\II{{\mathbf I}}
\def\JJ{{\mathbf J}}
\def\KK{{\mathbf K}}
\def\LL{{\mathbf L}}
\def\MM{{\mathbf M}}
\def\NN{{\mathbf N}}
\def\OO{{\mathbf O}}
\def\PP{{\mathbf P}}
\def\QQ{{\mathbf Q}}
\def\RR{{\mathbf R}}
\def\SS{{\mathbf S}}
\def\TT{{\mathbf T}}
\def\UU{{\mathbf U}}
\def\VV{{\mathbf V}}
\def\WW{{\mathbf W}}
\def\XX{{\mathbf X}}
\def\YY{{\mathbf Y}}
\def\ZZ{{\mathbf Z}}

\def\cA{{\mathcal A}}
\def\cB{{\mathcal B}}
\def\cC{{\mathcal C}}
\def\cD{{\mathcal D}}
\def\cE{{\mathcal E}}
\def\cF{{\mathcal F}}
\def\cG{{\mathcal G}}
\def\cH{{\mathcal H}}
\def\cI{{\mathcal I}}
\def\cJ{{\mathcal J}}
\def\cK{{\mathcal K}}
\def\cL{{\mathcal L}}
\def\cM{{\mathcal M}}
\def\cN{{\mathcal N}}
\def\cO{{\cmathal O}}
\def\cP{{\cmathal P}}
\def\cQ{{\mathcal Q}}
\def\cR{{\cmathal R}}
\def\cS{{\mathcal S}}
\def\cT{{\mathcal T}}
\def\cU{{\mathcal U}}
\def\cV{{\mathcal V}}
\def\cW{{\mathmathcal W}}
\def\cX{{\mathcal X}}
\def\cY{{\mathcal Y}}
\def\cZ{{\mathcal Z}}

\mathchardef\alphag="7C0B
\mathchardef\betag="7C0C
\mathchardef\gammag="7C0D
\mathchardef\deltag="7C0E
\mathchardef\varepsilong="7C22
\mathchardef\varphig="7C27
\mathchardef\psig="7C20
\mathchardef\zetag="7C10
\mathchardef\epsilong="7C0F
\mathchardef\rhog="7C1A
\mathchardef\taug="7C1C
\mathchardef\upsilong="7C1D
\mathchardef\iotag="7C13
\mathchardef\thetag="7C12
\mathchardef\pig="7C19
\mathchardef\sigmag="7C1B
\mathchardef\etag="7C11
\mathchardef\omegag="7C21
\mathchardef\kappag="7C14
\mathchardef\lambdag="7C15
\mathchardef\mug="7C16
\mathchardef\xig="7C18
\mathchardef\chig="7C1F
\mathchardef\nug="7C17
\mathchardef\varthetag="7C23
\mathchardef\varpig="7C24
\mathchardef\varrhog="7C25
\mathchardef\varsigmag="7C26
\mathchardef\Omegag="7C0A
\mathchardef\Thetag="7C02
\mathchardef\Sigmag="7C06
\mathchardef\Deltag="7C01
\mathchardef\Phig="7C08
\mathchardef\Gammag="7C00
\mathchardef\Psig="7C09
\mathchardef\Lambdag="7C03
\mathchardef\Xig="7C04
\mathchardef\Pig="7C05
\mathchardef\Upsilong="7C07

\def\barql{\bar \QQ_{\ell}}
\def\fd{F \!\hbox{-} \Delta}
\def\fdphi{F_{\varphi} \hbox{-} \Delta}
\def\ker{{\rm Ker} \,}
\def\st{{\rm St}}
\def\repr{{\rm Repr}}
\def\ind{{\rm Ind}}
\def\ch{{\rm ch} \,}
\def\chcl{{\rm chcl} \,}
\def\coker{{\rm Coker} \,}
\title[CHARACTER sums associated to finite Coxeter groups]{CHARACTER sums associated to finite Coxeter groups}
\author{Jan Denef}
\address{University of Leuven, Department of Mathematics,
Celestijnenlaan 200B, 3001 Leuven, Belgium}
\email{Jan.Denef@wis.kuleuven.ac.be}

\author{Fran\c cois Loeser}

\address{Centre de Math\'ematiques,
Ecole Polytechnique,
F-91128 Palaiseau
(URA 169 du CNRS), {\rm and}
Institut de Math\'{e}matiques,
Universit\'{e} P. et M. Curie, Case 247,
4 place Jussieu,
F-75252 Paris Cedex 05
(UMR 9994 du CNRS)}
\email{loeser@math.polytechnique.fr}

\date{}
\begin{abstract}The main result of this paper
is a character sum identity for Coxeter arrangements over finite fields
which is an analogue of Macdonald's conjecture
\cite{M} proved by Opdam \cite{O}.
\end{abstract}
\maketitle

\bigskip

\bigskip

%\sommaire
\setcounter{section}{-1}
\renewcommand{\theequation}{\thesection.\arabic{equation}}
\section{Introduction}Throughout this paper $\FF $
will denote a finite field of characteristic $p$
different from $2$.

We first slightly reformulate Macdonald's conjecture
in a form which has a direct analogue
over $\FF$.
Let $G$ be a finite subgroup of
${\rm GL}_{n} (\RR)$ generated by reflections and let $q$ be a
positive definite quadratic form which is invariant under $G$.
Let $\cA_{G}$ be the associated arrangement consisting of the reflection
hyperplanes.
Let $\ell_{1}, \ldots, \ell_{N}$ be equations for the
$N$ different
reflection hyperplanes.
Set $\Delta (x)= (\prod_{i = 1}^{N} \ell_{i})^{2}$
and let $d_{1}, \ldots, d_{n}$ be the degrees of $G$.
Macdonald's conjecture
\cite{M} is the following equality,
$$
\int_{\RR^{n}} \Delta (x)^{s} e^{-(\sum_{i}x_{i}^{2}) / 2}
dx =
(2 \pi)^{n / 2} \prod_{i= 1}^{n} \frac{\Gamma (d_{i} s + 1)}{\Gamma (s + 1)},
$$
when $q (x) = \sum_{i}x_{i}^{2}$ and the $\ell_{i}$ are normalized
in such a way
that
$\|\ell_{i}\| = \sqrt{2}$.
If we drop the condition that the $\ell_{i}$ are normalized, and if we
replace
$x$ by
$\sqrt{2} x$, we obtain
$$
\int_{\RR^{n}} \Delta (x)^{s} e^{-\sum_{i}x_{i}^{2}}
dx =
\pi^{n / 2} \left(\prod_{i = 1}^{N} \frac{{\|\ell_{i}\|}^{2}}{4}\right)^{s}
\prod_{i= 1}^{n} \frac{\Gamma (d_{i} s + 1)}{\Gamma (s + 1)}.
$$
Hence for $q$ an arbitrary
positive definite quadratic form which is invariant under $G$
we obtain
\begin{equation}
\int_{\RR^{n}} \Delta (x)^{s} e^{- q (x)}
dx =
\pi^{n / 2} \kappa^{s} \left(\prod_{i= 1}^{n} \frac{\Gamma (d_{i} s + 1)}
{\Gamma (s + 1)}\right)
({\rm discr} q)^{- \frac{1}{2}},
\end{equation}
with
\begin{equation}
\kappa = \prod_{i=1}^{N} \frac{q (\ell_{i})}{4}\, , 
\end{equation}
where  we consider $\ell_{i}$ in $q (\ell_{i})$
as a vector in $\RR^{n}$,
identifying
$\RR^{n}$ with its dual, by means of the quadratic form $q$.
Observe that
\begin{equation}
\sqrt{q (\ell_{i})} = {\rm Max}_{\sur{ x\in \RR^{n}}{q(x) = 1}}
\ell_{i}(x) \, . 
\end{equation}
Formula (0.1) can be reformulated as
\begin{equation}
\int_{K} \Delta (x)^{s} e^{- q (x)}
\sqrt{\Delta (x)}dx =
\pi^{n / 2} \kappa^{s + \frac{1}{2}}
\left(\prod_{i= 1}^{n} \frac{\Gamma (d_{i} (s + \frac{1}{2}))}{\Gamma (s
+
\frac{1}{2})}\right) \,
({\rm discr} q)^{- \frac{1}{2}},
\end{equation}
where $K = \RR^{n} / G$. Observe that
$\sqrt{\Delta (x)} dx$ is a $G$-invariant differential form.

Over $\FF$, there is the notion of a
Coxeter arrangement (1.4). It is
a triple,
$A = (V, G, q)$,
where $V$ is
a finite-dimensional $\FF$-vector space,
$G$ a finite subgroup of ${\rm GL} (V)$
generated by reflections,
and $q$ a $G$-invariant nondegenerate
symmetric bilinear form on $V$.
If $p$ does not divide $\vert G \vert$ one may define the degrees
of $G$, 
$d_{1}, \ldots, d_{n}$ (1.5). Let $\ell_{1}, \ldots, \ell_{N}$
be equations for the $N$ different
reflection hyperplanes.
Put $\Delta (x)= (\prod_{i = 1}^{N} \ell_{i})^{2}$.
Because
$p \not= 2$, we may define an element $\kappa$
of $\FF$ by (0.2).
Fix a non trivial additive character
$\psi : \FF \rightarrow \CC$.
The analogue of the integral in (0.4) will be the character sum
$$
S_{G} (\chi) := \sum_{x \in (U / G) (\FF)}
\chi (\Delta (x)) \psi (q (x))
$$
where
$\chi : \FF^{\times} \rightarrow \CC^{\times}$ is
a multiplicative character and $U$ denotes the complement of the
hypersurface
$\Delta (x) = 0$ in $V$. (We write
$q (x)$ for the quadratic form $q (x, x)$
associated to the  bilinear form $q$, and use the standard notation
$ (U / G)(\FF)$ to denote the set of $\FF$-rational points on the
quotient space $U / G$.)
The analogue of the Gamma function will be the Gauss sum $-g (\chi)$, where
$g (\chi) := -\sum_{x \in \FF^{\times}} \chi (x) \psi (x)$.
Our main result is the following.

\newtheorem*{MT}{Main Theorem}
\begin{MT}Let $A = (V, G, q)$ be a 
Coxeter arrangement over $\FF$. Assume that
$p$ does not divide
$\vert G \vert$.  Then
$\kappa \not=0$ and,
for every multiplicative
character $\chi : \FF^{\times} \rightarrow \CC^{\times}$,
$$
S_{G} (\chi) =
(- 1)^{n} \phi ({\rm discr} q)
g (\phi)^{n} \phi (\kappa) \chi (\kappa) \prod_{i = 1}^{n}
\frac{g ((\phi \chi)^{d_{i}})}{g (\phi \chi)} \, ,
$$
where $\phi$ denotes the unique  multiplicative character 
of $\FF$ of
order
2.
\end{MT}

Note the analogy with Macdonald's formula (0.4),
replacing $\sqrt \pi = \Gamma (1 / 2)$ by $- g (\phi)$.
(For some more remarks on this analogy, see (4.6) below.)
Our proof of the Main Theorem depends on 
Macdonald's formula (0.1) at one single point, namely for the proof of
Theorem 3.3 we need to know that,
with the notation of (0.1),
\begin{equation}
{\rm Max} \, \{ \Delta (x) \bigm | x \in \RR^{n}, q (x) = 1 \}
=
\kappa N^{-N} \prod_{i = 1}^{n} d_{i}^{d_{i}},
\end{equation}
which follows directly from (0.1), by elementary calculus.
Without relying on Macdonald's formula we can nevertheless
prove the Main Theorem for some $\kappa$ in $\FF \setminus \{0\}$
which is independent of $\chi$ without the assertion that $\kappa$ is
given by (0.2).

In the special case where $G$ is the symmetric group
$\cS_{n}$, the Main Theorem has been proved for all $p \not= 2$
by R. Evans \cite{E2}, extending the results of G. Anderson on Selberg
sums \cite{A1}. For the history and applications of such sums, we
refer to \cite{E1} and \cite{E3}. For recent work on character sums 
related to relative invariants of prehomogeneous vector spaces, see \cite{D-G}.

Our proof of the Main Theorem is entirely based on the cohomological
interpretation of character sums, using the Grothendieck-Lefschetz
trace formula (see \cite{SGA41/2} Sommes trig.
and \cite{K}).
In the present situation, the cohomology is concentrated in middle
dimension
(Proposition 2.1.2) and has rank 1 (Corollary 4.3.2
and formula 4.4.6),
so that we have only to calculate the determinant
of the Frobenius action on the cohomology, which is done by Laumon's
product formula \cite{La} 3.2.1.1.

\bigskip
The plan of the paper is the following. In the first section
we study Coxeter arrangements over finite fields and their liftings to
characteristic zero.
In section 2 we study character sums associated to central arrangements
that are invariant under a finite group. This will allow us to evaluate
$S_{G} (\phi)$. In the next section we review  results
on
monodromy and critical values for Coxeter arrangements from
\cite{D-L3} and \cite{D-L4}.
In the last section, using all the previous results, we prove the Main Theorem.

\bigskip
{\small Most of the present  work dates from 1992. It has been presented,
under some additional hypotheses, in lectures by the first named
author in Orsay (SAGA) in November 1992 and in Amsterdam in January 1993.
}

\renewcommand{\theequation}{\arabic{equation}}
\setcounter{equation}{0}
\section{Liftable Coxeter arrangements over $\FF$}
\subsection{}Let $V$ be a finite dimensional vector space over a
field. By a {\em hyperplane
arrangement} in $V$
we mean a finite set of affine
hyperplanes in $V$.
If all the hyperplanes contain $0$, the arrangement is said to
be {\em central}. We call an endomorphism of $V$ a {\em reflection}
if it has order 2 and fixes pointwise some hyperplane.

We define a {\em classical Coxeter arrangement}
as a triple
$A = (V, G, q)$,
where $V$ is a finite dimensional vector space over $\RR$,
$G$ is a finite subgroup of ${\rm GL} (V)$
generated by reflections,
and $q$ is a $G$-invariant {\em positive definite}
symmetric bilinear form on $V$.
We define a {\em Coxeter arrangement over $\CC$}
as a triple
$A = (V, G, q)$,
where $V$ is a finite dimensional vector space over $\CC$,
$G$ is a finite subgroup of ${\rm GL} (V)$
generated by reflections,
and $q$ is a $G$-invariant nondegenerate
symmetric bilinear form on $V$,
which arises by extension of scalars from a classical
Coxeter arrangement.
(In fact, by
the argument given at the end of the proof of Proposition 1.6,
the last condition in this definition is automatically
verified.)

Similarly, when $T$ is a subring of $\CC$ with fraction field $K$, a
{\em
Coxeter
arrangement over $T$} will be a triple
$A = (M, G, q)$
where $M$ is
a free $T$-module of finite rank,
$G$ a finite subgroup of ${\rm GL} (M)$
generated by reflections,
and $q$ a $G$-invariant nondegenerate
symmetric bilinear form on $M$,
which induces  by extension of scalars  a Coxeter arrangement over $\CC$.

The finite set of reflection hyperplanes in $M \otimes_{T} K$
defines a central
hyperplane arrangement in $M \otimes_{T} K$
which we denote by ${\cal A}_G$.

\subsection{}Let $T$ be a discrete valuation
ring with fraction field
$K$. Let $\got P$ be the maximal ideal of $T$
with residue field
$\overline K = T / \got P$. Let $M$ be a free
$T$-module of finite rank, and set
$V = M \otimes_{T} K$, $\overline V = M \otimes_{T} \overline K$.
We denote 
the reduction map $M \rightarrow \overline V$
by $x \mapsto \overline x$.
For any linear subspace $W$ in $V$, we
denote by $\overline W$ the reduction modulo $\got P$
of $W$ which  is defined by
$$
\overline W = \{\overline x \bigm| x \in W \cap M\}.
$$
Observe that ${\rm dim}_{K} W = {\rm dim}_{\overline K} \overline W$.

\begin{definition}Let $\cA = (H_{i})_{i \in J}$ be a central
arrangement in $V$. We say
that $\cA$ has good reduction mod $\got P$
if, for any $I \subset J$,
$$
\overline{\cap_{i \in I}H_{i}}
=
\cap_{i \in I}\overline{H_{i}}.
$$
\end{definition}

For $\cA$ an arrangement in  $V$
consider the arrangement $\overline \cA$
consisting of the hyperplanes $\overline H$ with $H$ in $\cA$.
Denote 
the lattices of $\cA$ and $\overline \cA$ by $L (\cA)$ and $L (\overline \cA)$.
There is a canonical inclusion preserving map
$$
\theta : \left\{ \begin{array}{l}
L (\cA) \longrightarrow L (\overline \cA)\\
X \mapsto \cap_{\sur{H \in \cA}{H \supset X}} \overline H.
\end{array}\right.
$$

\begin{prop}If $\cA$ has good reduction mod $\got P$,
then $\theta$ is an isomorphism of lattices.
\end{prop}
\begin{proof} This is straightforward because $\theta (X) = \overline X$
when $\cA$ has good reduction. Also the equality ${\rm dim}_{K} W = {\rm
dim}_{\overline K} \overline W$ holds
for any linear subspace $W$ of $V$.\end{proof}

\subsection{}We assume the notation of (1.2).
Suppose that $T \subset \CC$ and let
$(M, G, q)$ be a Coxeter arrangement over $T$.
Thus, in particular, $G \subset {\rm GL} (M)$.
Let $\cA_{G}$ be the central arrangement in $V$ consisting of all the
reflection hyperplanes of $G$.
We denote  the image of $G$ in
${\rm GL} (\overline V)$ by
$\overline G$ and we choose linear forms $\ell_{i}$
over $T$ which define the hyperplanes of 
$\cA_{G}$ and which are not zero modulo $\got P$.
Let $\kappa$ be as in (0.2).

\begin{prop}Assume that
$\got P$ does  not contain $\vert G \vert$. Then $\cA_{G}$ has good
reduction mod
$\got P$. Furthermore the canonical morphism
$G \rightarrow \overline G$
is an isomorphism,
and $G$ acts freely on the complement
in
$\overline V$ of the hyperplanes in $\overline {\cA_{G}}$.
\end{prop}

\begin{proof}For the first statement it suffices to show that for all $I \subset J$
$$
{\rm dim} \, \cap_{i \in I}\overline{H_{i}}
\leq
{\rm dim} \, \overline{\cap_{i \in I}H_{i}}.
$$
We may suppose that the $H_{i}$'s with $i$ in $I$
are linearly independent. Let $s_{i}$ be the reflection
with respect to $H_{i}$
and set $\gamma = \prod_{i \in I} s_{i}$,
${\rm Fix} (\gamma) =
\{ x \in V \bigm| \gamma (x) = x \}$ and
${\rm Fix} (\overline \gamma) =
\{ \overline x \in \overline V \bigm| \overline \gamma (\overline x) =
\overline x \}$.
From the proof of Theorem 6.27 (2) on p.225 of \cite{O-T}
we see that ${\rm Fix} (\gamma) = \cap_{i \in I}H_{i}$.
Since $\overline H_{i} \subset {\rm Fix} (\overline s_{i})$, 
${\rm dim} \, \cap_{i \in I}\overline{H_{i}} \leq {\rm dim} \, {\rm Fix}
(\overline \gamma)$.
But this last dimension is equal to the number of eigenvalues of
$\overline \gamma$ which are equal to 1, {\it i.e.} the number of eigenvalues of
$\gamma$ which are congruent  to 1 mod $\got P$
(counting multiplicities and enlarging $T$ so that it contains all
eigenvalues of $\gamma$).
But if $\lambda$ is an eigenvalue of
$\gamma$ which is
congruent  to 1 mod $\got P$, then $\lambda = 1$,
because $\lambda^{\vert G \vert} = 1$, and
$\got P$ does  not divide $\vert G \vert$.
But the number of eigenvalues of
$\gamma$ which are equal to 1
(counting multiplicities) is equal to ${\rm dim} \, {\rm Fix}
(\gamma)$.
We now deduce the first statement from the equalities
$$
{\rm dim} \, {\rm Fix}
(\gamma)
=
{\rm dim} \, \cap_{i \in I}{H_{i}}
=
{\rm dim} \, \overline {\cap_{i \in I}{H_{i}}}.
$$
For the second statement we just have to observe
that if an element $\gamma$ of $G$
is in the kernel of the morphism
$G \rightarrow \overline G$,
all the eigenvalues of
$\gamma$ are congruent  to 1 mod $\got P$, hence are equal to 1
by the above argument. The order of $\gamma$ being finite,
we conclude that $\gamma$ is equal to 1.
Finally the third assertion follows from the well known fact
({\it cf.} \cite{B} Ch.V, \S\kern .15em 3 Proposition 3)
that $G$ acts freely on the complement in $V$
of the hyperplanes in $\cA_{G}$, and from the observation that
${\rm Fix} (\overline \alpha)
= \overline {{\rm Fix} (\alpha)}$, for all $\alpha$ in $G$.\end{proof}

\medskip
For each reflection $s \in G$ we choose an
eigenvector $\alpha_{s} \in M$ of $s$
which is orthogonal
(with respect to $q$) to the reflection hyperplane of $s$.

Set
$\Phi = \{\alpha_{s} \, \vert \, s \, \hbox{reflec\-tion in} \,  G \}$. Note that
\begin{equation}
s (x) = x - 2\frac{q (x, \alpha_{s})}
{q (\alpha_{s}, \alpha_{s})} \alpha_{s} \quad \hbox{for all} \, x \in V.
\end{equation}
If $q (\alpha, \alpha) \not \in \got P$
for each $\alpha \in \Phi$, then,
\begin{equation}
\overline s (\overline x) = \overline x - 2\frac{\overline q (\overline x,
\overline \alpha_{s})}
{\overline q (\overline \alpha_{s}, \overline \alpha_{s})} \overline
\alpha_{s} \quad
\hbox{for all} \, \overline x \in \overline V.
\end{equation}

Denote by $\overline q$ the symmetric bilinear form induced by
$q$ on $\overline V$.
\begin{lem}Assume that $G$ is essential (meaning that
$V^{G} = \{0\}$), that
$\got P$ does  not contain $\vert G \vert$,
and that $q (\alpha, \alpha) \not \in \got P$
for all $\alpha \in \Phi$. Then the bilinear form $\overline q$ is
nondegenerate.
\end{lem}

\begin{proof}Assume that there exists $\overline x \in \overline V$
with $\overline q (\overline x, \overline y) = \overline 0$ for all
$\overline y \in \overline V$.
Then (2) implies that $\overline s (\overline x) = \overline x$ for every
reflection $s$ in $G$.
Let $H_{s}$ be the reflection hyperplane of the reflection $s$.
Then $\overline x \in \overline H_{s}$ because otherwise
$\overline s$ is the identity and all eigenvalues of $s$ would be congruent to
1 mod $\got P$. Hence $\overline x \in \cap_{H \in \cA_{G}} \overline H
= \overline{\cap_{H \in \cA_{G}} H} = \overline{V^{G}} = \overline 0$
by Proposition 1.3.1.\end{proof}

\begin{prop}Assume that
$\got P$ does  not contain $\vert G \vert$
and that the bilinear form $\overline q$ is nondegenerate. Then
$\kappa \in T \setminus \got P$, and, for any
$I \subset J$, the restriction of $\overline q$ to
$\cap_{i \in I} \overline H_{i}$ is nondegenerate.
\end{prop}

\begin{proof}Observe first that the characteristic of $\overline K$ is not
2 because
$2 \, \vert \, \vert G \vert$. For each reflection $s \in G$, choose an
eigenvector
$\alpha_{s} \in M$ of $s$ which is orthogonal (with respect to $q$) to the
reflection
hyperplane of $s$ such that $\overline \alpha_{s} \not = 0$.
Since $\overline q$ is nondegenerate there exists $x \in M$
such that $\overline q (\overline x, \overline \alpha_{s}) \not =  0$.
Thus
$q (x, \alpha_{s}) \not \equiv 0 \, \hbox{mod} \, \got P$.
Because $s (x) \in M$, formula (1) yields $q (\alpha_{s},
\alpha_{s}) \not \equiv 0 \, \hbox{mod} \, \got P$.
This implies that $\kappa \in T - \{0\}$.
Set $U = \cap_{i \in I} H_{i}$, thus
$\overline U = \cap_{i \in I} \overline H_{i}$.
We have to show that the restriction of $\overline q$ to $\overline
U$ is nondegenerate. Because
$\overline q$ is nondegenerate, it suffices to prove that the
restriction to
$(\overline U)^{\perp}$ is nondegenerate.
We have
$\overline {(U^{\perp})} = (\overline
U)^{\perp}$ because,
clearly,
$\overline {(U^{\perp})} \subset (\overline U)^{\perp}$
and both have the same dimension $n - \hbox{dim} \, U$.
Thus it remains to prove that
$\overline q$
is nondegenerate on $\overline {(U^{\perp})}$.
Let $W$ be the subgroup of $G$ which is generated by the reflections
with respect to the reflection 
hyperplanes
of $G$ which contain $U$.
Clearly, $W$ is a finite Coxeter group in $U^{\perp}$ and $W \subset
{\rm GL} (U^{\perp} \cap M)$.
It is easy to see that the hypothesis of Lemma 1.3.2 is satisfied if we
replace $V$ by
$U^{\perp}$, $M$ by $U^{\perp} \cap M$ and $q$ by its restriction to
$U^{\perp} \cap M$. Indeed
$(U^{\perp})^{W} = \{ 0\}$ because $V^{W}$ is contained in $U$, and thus
$(U^{\perp})^{W} \subset U^{\perp} \cap U = \{ 0\}$ since $q_{\vert U}$ 
arises from a positive definite form over $\RR$.
Hence Lemma 1.3.2 shows that $\overline
{(q_{\vert U^{\perp} \cap M})}$ is nondegenerate. Thus
$\overline q_{\vert \overline{(U^{\perp})}}$ is also
nondegenerate.\end{proof}

\subsection{}A {\em Coxeter arrangement over $\FF$}
is a triple
$A = (V, G, q)$
where $V$ is
a finite dimensional $\FF$-vector space,
$G$ a finite subgroup of ${\rm GL} (V)$
generated by reflections,
and $q$ a $G$-invariant nondegenerate
symmetric bilinear form on $V$.
By a reflection we mean an endomorphism of $V$
of order 2, which fixes pointwise some hyperplane.
The finite set of reflection hyperplanes
defines a central
hyperplane arrangement
which we denote by ${\cal A}_G$.
A Coxeter arrangement $A = (V, G, q)$
over $\FF$ will be called {\em liftable}
if there exists an embedding $T \hookrightarrow \CC$
of a 
discrete valuation ring $T$ with 
residue field $\FF$ and fraction field of finite degree
over $\QQ$,
and a Coxeter arrangement $A_{T} = (M_T, G_T, q_T)$
over $T$
such that $V$ is isomorphic to
$M_T \otimes_T \FF$ and, after identification of
$M_T \otimes_T \FF$ with $V$,
$G = \overline G_T$ and $q = \overline q_T$. Here we denote 
the symmetric bilinear form
induced by $q_T$ on $V$ by $\overline q_T$.
Such a Coxeter arrangement $A_{T}$ will be called a lifting of
$A$ over $T$.

The following statement is a direct consequence of Proposition 1.3.1.
\setcounter{stheorem}{3}
\begin{scor}Let $T$ be a 
discrete valuation ring
with 
residue field $\FF$ 
and fraction field of finite degree
over $\QQ$. Fix an embedding $T \hookrightarrow \CC$.
Let $A = (M, G, q)$ be a Coxeter arrangement
over $T$. Suppose that
$p$ does not divide
$\vert G \vert$
and that $\overline q$, the symmetric bilinear form
induced
by $q$ on $\overline V$, is nondegenerate.
Then $\overline A := (\overline V, \overline G, \overline q)$
is a liftable Coxeter arrangement
over $\FF$.
Furthermore the arrangement ${\cal A}_G$ has good reduction and we have
$\overline {{\cal A}_G} = {\cal A}_{\overline G}$.
\hfill $\qed$
\end{scor}

\subsection{}Let $A = (V, G, q)$ be a Coxeter arrangement
over
$\FF$. Let us denote by $S (V)$ the symmetric algebra
of $V$ and by $R(V)$ the subalgebra of $G$-invariants.
If $p$ does not divide $\vert G \vert$, by Chevalley's Theorem
(\cite{B} Ch.V, \S\kern .15em 5 Th\'{e}or\`{e}me 4),
$R (V)$ is 
generated as an $\FF$-algebra by $n := {\rm dim} \, V$
algebraically independent homogeneous polynomials (together with 1).
The degrees of these polynomials are called
the degrees
$d_{1}, \ldots, d_{n}$ of $A$ (or the degrees of $G$).
Their product being equal to
$\vert G \vert$, they are also prime to $p$ ({\it loc. cit.} p.115).

\begin{sprop}Let $A = (V, G, q)$ be a liftable
Coxeter arrangement
over
$\FF$. Assume that $p$ does not divide $\vert G \vert$.
The degrees of $A$ coincide with those of any
lifting of $A$ viewed as a complex
Coxeter arrangement.
\end{sprop}

\begin{proof}Let $T \subset \CC$ be a 
discrete valuation ring
having $\FF$
as residue field and let $A_{T}$ be a lifting of
$A$ over $T$. We have to prove that the degrees of
$A$ and $A_{T}$ (viewed as a complex
Coxeter arrangement)
are the same. For that we may clearly suppose that $T$ is complete.
Let $p_{1}, \ldots, p_{n}$ be homogeneous generators of
$R (V)$. As $p$ does not divide $\vert G \vert$,
$p_{1}, \ldots, p_{n}$ can clearly be lifted, by averaging,
to elements
$q_{1}, \ldots, q_{n}$ of
$R (M_{T})$, the subalgebra of $G$-invariants of the symmetric algebra
of $M_{T}$. Let $\pi$ be a uniformizing parameter of
$T$. The polynomials $q_{1}, \ldots, q_{n}$ generate the algebra
$R (M_{T})$
modulo $\pi$, hence, by induction, modulo $\pi^{i}$ for any $i$.
Passing to the limit, we see  that
$q_{1}, \ldots, q_{n}$ generate the algebra
$R (M_{T})$, which implies the result.\end{proof}

\subsection{}By the following proposition, we have liftability
of Coxeter arrangements over finite fields, as soon as the characteristic
does not divide the order of the Coxeter group.

\begin{sprop}Let $A = (V, G, q)$ be a 
Coxeter arrangement
over
$\FF$. If $p$ does not divide $\vert G \vert$, then $A$ is liftable.
\end{sprop}

\begin{proof}Since the group ring $\FF [G]$ is semisimple, each simple
component of the $\FF [G]$-module $V$ is a direct summand of
$\FF [G]$. Thus $V$ is a projective $\FF [G]$-module,
and can hence be lifted to a $T [G]$-module $M$, of finite rank over
$T$, for some discrete valuation ring $T \subset \CC$, 
with residue field $\FF$ and fraction field $K$
of finite degree over $\QQ$.
In fact, by
\cite{Se}  \S\kern .15em 14.4
Proposition 42, we can lift $V$
to a $R [G]$-module $M$, of finite rank over
$R$, for some complete discrete valuation ring $R \subset \CC$, 
with residue field $\FF$, and it follows from
a result of M. Greenberg \cite{G}
that we can work with a henselian discrete
valuation ring instead of a complete one as in {\it loc. cit}.
[As was suggested by the referee one could also
observe that the group algebra
$\ZZ [\frac{1}{|G|}][G]$ is an Azumaya
algebra over its center which is unramified over $ \ZZ [\frac{1}{|G|}]$,
and take $T$ to split them ({\it cf.} \cite{Mi} Ch.IV \S\kern .15em 1).]
Since the reduction modulo $p$
is injective on the roots of unity of order dividing
$|G|$, we see that the reflections
in $G \subset {\rm GL} (V)$ can be lifted to reflections
in ${\rm GL} (M)$. By averaging, the form $q$ can be lifted to
a $G$-invariant nondegenerate symmetric bilinear form $q_{T}$
on $M$. Moreover, the action of $G$ on $M \otimes_{T} \CC$
can be realized over $\RR$
(see, {\it e.g.},  \cite{Se} \S\kern .15em 13.2
Th\'{e}or\`{e}me 31) in such a way that $q_{T}$ corresponds to a
positive definite symmetric bilinear form over $\RR$
(indeed, this follows from the proof in {\it loc. cit.}).
Thus $(M, G, q_{T})$ is a Coxeter arrangement over $T$.\end{proof}

\section{Character sums associated to invariant central arrangements}
\setcounter{subsection}{-1}
\setcounter{equation}{0}
\subsection{The data}Let $p$ be a prime number different
from 2 and let
$\FF$ be a finite field of characteristic $p$. We denote by $\bar \FF$
a fixed algebraic closure of $\FF$.
Let $f = \prod_{i \in I} \ell_{i}$ be a product of (not necessarly distinct)
linear forms $\ell_{i}$ on $\AA^{n}_{\FF}$.
The locus of $f = 0$ is an arrangement $\cA$. We denote by
$U$ its complement in $\AA^{n}_{\FF}$.
Let $q$ be a quadratic form on $\AA^{n}_{\FF}$ whose restriction to any
stratum of $\cA$ is nondegenerate.
Choose a prime number $\ell \not= p$.
Let $\chi$ be a multiplicative character
$\FF^{\times} \rightarrow \barql^{\times}$
and let
$\psi$ be a non trivial additive character
$\FF \rightarrow \barql^{\times}$. (Here
$\barql$ denotes an algebraic closure of the field
of $\ell$-adic numbers.)
We denote by $\cL_{\chi}$ and $\cL_{\psi}$
the corresponding Kummer and Artin-Schreier sheaves
(see \cite{SGA41/2} Sommes trig.).
Let $G$ be a finite group, leaving
$f$ and
$q$ invariant, and acting freely on $U$.
We consider the exponential sum
\begin{align*}
S_{G} (\chi) & =  \sum_{x \in (U / G) (\FF)}
\chi (f(x)) \psi (q(x))\\
& = 
{\rm Tr} \, (F,
(H^{\cdot}_{c} (U_{\bar \FF},
f^{\ast} \cL_{\chi} \otimes q^{\ast} \cL_{\psi}))^{G}).
\end{align*}
Here $F$ denotes the geometric Frobenius automorphism
and the second equality follows from Grothendieck's trace formula,
see {\it loc. cit}.

At this point we should perhaps recall how
the action of $G$
on
$H^{\cdot}_{c} (U_{\bar \FF},
f^{\ast} \cL_{\chi} \otimes q^{\ast} \cL_{\psi})$ is defined.
It is induced
by the canonical isomorphism
$$
H^{\cdot}_{c} (U_{\bar \FF},
f^{\ast} \cL_{\chi} \otimes q^{\ast} \cL_{\psi})
\simeq
H^{\cdot}_{c} ((U / G)_{\bar \FF},
\pi_{\ast} \barql \otimes f^{\ast}_{G} \cL_{\chi}
\otimes q^{\ast}_{G} \cL_{\psi})
$$
where $\pi : U \rightarrow U / G$ is the natural map,
$f_{G}$ and $q_{G} :U / G \rightarrow \AA^{1}$ are
induced by
$f$ and $q$, with the natural action of $G$
on $\pi_{\ast} \barql$
and the trivial action on
$f^{\ast}_{G} \cL_{\chi}$ and
$q^{\ast}_{G} \cL_{\psi}$.

\medskip

\subsection{Concentration of the cohomology in the middle
dimension}We construct the following compactification $\tilde q$
of $q$.
Let $\Gamma$ be the closure in
$\PP^{n} \times \AA^{1}$ of the graph of
$q$ in $\AA^{n} \times \AA^{1}$ and
let
$\tilde q$ be the projection
$\Gamma \subset \PP^{n} \times \AA^{1} \rightarrow \AA^{1}$.
Clearly $\tilde q$ is proper.
We denote by $j : U \rightarrow \Gamma$ the canonical open immersion.
Let $\cF = j_{!} f^{\ast} \cL_{\chi}$.

\begin{lem}Outside
$0 \in \AA^{n} \subset \Gamma$,
the morphism
$\tilde q$ is locally acyclic\footnote{For the notion of ``locally acyclic",
see \cite{SGA41/2} Th. finitude 2.12.}
with respect to the sheaf $\cF$
and $\Gamma$ is smooth over $\FF$.
\end{lem}

\begin{proof}We look locally at $a \in \Gamma \setminus \{0\}$.

\noindent Case (i): $a \in \AA^{n} \subset \Gamma$. We may
assume that the linear forms $\ell_{i}$ which vanish at $a$
are $x_{n}, x_{n-1}, \ldots, x_{r}$
and linear combinations of these, with $r > 1$.
The restriction of $q$ to $x_{n} = x_{n - 1} = \cdots
= x_{r} = 0$ is nondegenerate
and $a \not= 0$, thus the restriction of $q$ is nonsingular at $a$.
Hence $x_{n}, x_{n-1}, \ldots, x_{r}, q - q (a)$ are part of a system
of local parameters for $\Gamma$ at $a$. In the diagram
\begin{equation*}\xymatrix{
{\rm Spec} \, \FF [x_{r}, \ldots, x_{n}] \ar[d] &\ar[l]
{\rm Spec} \, \FF [x_{r}, \ldots, x_{n}, q] \ar[d]
& & \ar[ll]_-{\alpha = (x_{r}, \ldots, x_{n}, q)} \AA^{n} \\
{\rm Spec} \, \FF & \ar[l] {\rm Spec} \, \FF [q]
}
\end{equation*}
the square is cartesian, the map $\alpha$
is smooth at $a$ and, locally at $a$ for the etale topology,
the sheaf $\cF$ is the pull back of a sheaf on 
${\rm Spec} \, \FF [x_{r}, \ldots, x_{n}]$,
thus $q$
is locally acyclic for $\cF$ at $a$.

\noindent Case (ii):  $a \notin \AA^{n} \subset \Gamma$.
The equation for $\Gamma$ in $\PP^{n} \times \AA^{1}$
is $y x_{0}^{2} = q (x_{1}, \ldots, x_{n})$
where $x_{0}, x_{1}, \ldots, x_{n}$ are projective coordinates on
$\PP^{n}$ and $y$ is the coordinate on $\AA^{1}$.
The map $\tilde q$ is given by
$(x, y) \mapsto y$.
Let the coordinates of $a$ be
$x_{i} = \alpha_{i}$ and $y = b$.
Since $a \notin \AA^{n}$
we have $\alpha_{0} = 0$.
Hence $q (\alpha_{1}, \ldots, \alpha_{n}) = 0$.
We may assume that, {\it e.g.},  $\alpha_{1} \not= 0$.
We may assume also that the linear forms $\ell_{i}$ which vanish at
$(\alpha_{1},
\ldots, \alpha_{n})$
are $x_{n}, x_{n-1}, \ldots, x_{r}$
and linear combinations of these, with $r > 1$.
Since the restriction of $q$
to
$x_{n} = x_{n - 1} = \cdots
= x_{r} = 0$ is nondegenerate,
$(\alpha_{1},
\ldots, \alpha_{n})$
is a nonsingular projective point
of the projective quadric
$q (x_{1}, \ldots, x_{r - 1}, 0, \ldots, 0) = 0$ and $r > 2$.
Hence
$(1, \frac{\alpha_{2}}{\alpha_{1}}, \frac{\alpha_{3}}{\alpha_{1}},
\ldots)$ is a nonsingular point of the affine quadric
$q (1, t_{2}, \ldots, t_{r - 1}, 0, \ldots, 0) = 0$.
Moreover, $\Gamma$ is defined
in a neighbourhood of
$a$ by the affine equation
$y t_{0}^{2} = q (1, t_{2}, \ldots, t_{n})$.
Here $t_{0}$, $t_{2}$, $\ldots$, $t_{n}$ are coordinates on $\AA^{n}$
and $y$ is the coordinate on $\AA^{1}$.
We see that $t_{0}$, $y - b$,
$t_{r}$, $\ldots$, $t_{n}$ are part of a system of local parameters
for $\Gamma$ at $a$. Since $\Gamma \setminus U$ is an hypersurface
of $\Gamma$ given locally  at $a$ by an equation only involving
$t_{0}$, $t_{r}$, $\ldots$, $t_{n}$,
we may conclude as in case (i) that $\tilde q$
is locally acyclic for $\cF$ at $a$.\end{proof}

\begin{prop}Assume that the sheaves
$R^{i}q_{\vert U !} f^{\ast} \cL_{\chi}$ have tame ramification at $\infty$
for all integers $i$.
Then
\begin{equation*}
H^{i}_{c} (U_{\bar \FF},
f^{\ast} \cL_{\chi} \otimes q^{\ast} \cL_{\psi}) = 0 \quad
\hbox{if} \quad
i \not= n.
\end{equation*}
\end{prop}

\begin{proof}It is a direct consequence of Lemma 2.1.1
and of a straightforward
adaptation
of \cite{D-L1} Proposition 3.1. Indeed
part 3.1.1  of that proposition remains valid when we
replace the sheaf
$\QQ_{\ell}$ by any constructible
$\barql$-sheaf.
Thus by Lemma
2.1.1
we get
$H^{i}_{c} (U_{\bar \FF},
f^{\ast} \cL_{\chi} \otimes q^{\ast} \cL_{\psi})
=
H^{i}_{c} (\Gamma_{\bar \FF},
\cF \otimes \tilde q^{\ast} \cL_{\psi})
=
0$
when $i > n$ provided
that
the sheaves
$R^{i}q_{\vert U !} f^{\ast} \cL_{\chi}$ have tame ramification at $\infty$
for all $i$. The case $i <n$
follows by Poincar\'{e} duality
and the affineness of $U$.\end{proof}

\medskip

\subsection{The case when $\chi$
is trivial}We assume throughout
this paragraph that $\chi$ is trivial,
thus
$$
S_{G}   = \sum_{x \in (U / G) (\FF)}
\psi (q(x)).$$
We also consider
$$S = \sum_{x \in U (\FF)}
\psi (q(x)) =
{\rm Tr} \, (F,
H^{\cdot}_{c} (U_{\bar \FF},
q^{\ast} \cL_{\psi})).
$$
We may assume  that the $\ell_{i}$'s
define different hyperplanes.

\begin{lem}If $\cA = \emptyset$ and
$U = \AA^{n}$,
then
$$
S = \phi ({\rm discr} q) (- g (\phi))^{n}
\, ,
$$
where $\phi$ denotes the unique multiplicative character of
$\FF$ of
order
2.
\end{lem}

\begin{proof}Reduce $q$ to diagonal form and use
$\sum_{x \in \AA^{1} (\FF)} \psi (x^{2}) = - g (\phi)$.\end{proof}

\begin{lem}The following formula holds,
\begin{equation*}
H^{i}_{c} (U_{\bar \FF},
q^{\ast} \cL_{\psi}) = 0 \quad
\hbox{if} \quad
i \not= n.
\end{equation*}
\end{lem}

\begin{proof}For $U = \AA^n$, this is well known
(reduce $q$ to diagonal form). The general case
is a direct consequence of Lemma 3.8.1 of \cite{D-L1}, and of the fact that
$U$ is affine.\end{proof}

\begin{lem}The vector space
$H^{n}_{c} (\AA^{n}_{\bar \FF}, q^{\ast} \cL_{\psi})$
is one dimensional and the action
of an element
$\sigma$ of $G$ on this space is given by multiplication by
its determinant (as an operator
on
$\AA^{n}$), 
${\rm det} \, \sigma = \pm 1$.
\end{lem}

\begin{proof}The first statement is well known.
For the second statement,
we may assume $\sigma$ is a reflection since the orthogonal group is
generated by reflections. By  K\"{u}nneth, the situation is reduced to
the case where $n = 1$ and then it is easy.\end{proof}

\medskip

\noindent {\it Notation}. --- We set
$$
{\rm det} \, (U, G, q) :=
{\rm det} \, (F, H^{\cdot}_{c} (U_{\bar \FF},
q^{\ast} \cL_{\psi})^{G})
=
{\rm det} \, (F, H^{n}_{c} (U_{\bar \FF},
q^{\ast} \cL_{\psi})^{G})^{(- 1)^{n}},
$$
with $n = {\rm dim} \, U$,
and we set $\nu_{U, G} = 0$
if there is a $\sigma$ in $G$ with ${\rm det} \, \sigma \not = 1$,
and $\nu_{U, G} = 1$ otherwise, with ${\rm det} \, \sigma$
denoting the determinant of the action of $\sigma$ on the affine space
of $U$.
Let $L (\cA)$, also denoted by $L (U)$, be the lattice
of the arrangement
$\cA$. For $X \in L (\cA)$ we set $X^{\circ} = X \setminus
\cup_{\sur{Y \in L (\cA)}{X \not\subset Y \hfill}}Y$
and 
$\st_{G} (X) = \{\sigma \in G \bigm| \sigma (X) \subset X\}$.
We will use the above notation also when $G$ does not act freely on $U$,
for example when $U$ is replaced by
$X^{\circ}$ and $G$ by
$\st_{G} (X)$. From now on till the end of subsection (2.2), we do not assume
that $G$ acts freely on $U$.

\begin{theorem}The following formula holds,
$$
{\rm det} \, (U, G, q) =
(\phi ({\rm discr}\, q)
g (\phi)^n)^{(- 1)^{n} \nu_{U, G}}
\prod_{X \in L (U) / G}
{\rm det} \, (X^{\circ}, \st_{G} (X), q_{\vert X})^{-1}.
$$
\end{theorem}

\begin{proof}The theorem follows directly from Lemmas 2.2.1, 2.2.2 and 2.2.3
and the decomposition
$U = \AA^{n} \setminus
\coprod_{\bar X \in L (U) / G} \coprod_{X \in \bar X} X^{\circ}$.\end{proof}

\medskip

\noindent {\it Notation}. --- We denote by
$\ch [\AA^{n}, X]$, for $X \in L (\cA)$, the set of all chains in the
lattice
$L (\cA)$ which connect $\AA^{n}$ to $X$.
By definition, such a chain is a sequence
$X {\subset} X_1 {\subset}
X_2 {\subset} \cdots {\subset}
X_m {\subset}
\AA^n$, with $X_i$ in $L (\cA)$, the inclusions being strict.
Note that $\st_{G} (X)$ acts on
$\ch [\AA^{n}, X]$, so that we can consider
the quotient
$\chcl [X] := \ch [\AA^{n}, X] / \st_{G} (X)$.
For $c$ in $\chcl [X]$
we denote by $|c|$ the length of a chain representing $c$.
For $X$ in $L (\cA)$, we let $\repr_{X}$
be the following virtual representation of
$\st_{G} (X)$
$$
\repr_{X} :=
(- 1)^{n - {\rm dim}\, X}
\sum_{c \in \chcl [X]} (- 1)^{|c| - 1}
\ind_{\st_{G} (c)}^{\st_{G} (X)} ({\rm det}_{\vert X}),
$$
where
$\ind_{\st_{G} (c)}^{\st_{G} (X)} ({\rm det}_{\vert X})$
is the induced representation of the representation
which to an element in
$\st_{G} (c)$ assigns its determinant as an operator on $X$,
and $\st_{G} (c)$ is the stabilizer of a chosen representative of $c$
in $\ch [\AA^{n}, X]$.

\begin{theorem}There is an isomorphism of virtual
representations
of $G$,
$$
H^{n}_{c} (U_{\bar \FF}, q^{\ast} \cL_{\psi})
\simeq
\sum_{X \in L (\cA) / G}
\ind^{G}_{\st_{G} (X)} (\repr_{X}).
$$
Moreover this isomorphism
preserves the Frobenius action if we let it
act by multiplication by
$(\phi ({\rm discr}\, q_{\vert X})) g (\phi)^{{\rm dim} \, X}$
on
$\repr_{X}$.
\end{theorem}

\begin{proof}By using the same decomposition of $U$ as in the proof
of Theorem 2.2.4,
we deduce from Lemma 2.2.2
the following isomorphism of virtual representations of $G$,
$$
H^{n}_{c} (U_{\bar \FF}, q^{\ast} \cL_{\psi})
\simeq
\repr_{\AA^{n}}
-
(- 1)^{n}
\sum_{\sur{X \in L (\cA) / G}{X \not= \AA^{n} \hfill}}
\ind^{G}_{\st_{G} (X)} (-1)^{{\rm dim} \, X}
H^{{\rm dim} \, X}_{c} (X^{\circ}_{\bar \FF}, q^{\ast} \cL_{\psi}),
$$
which is compatible with the Frobenius action.
The theorem follows now from Lemmas
2.2.1 and 2.2.3, by recursion and transitivity of
induction.\end{proof}

\medskip
The remaining results in this subsection 2.2 will not be used in the rest of
the paper.

\medskip

\noindent {\it Notation}. --- For $X \in L (\cA)$,
we set
$d (X) := {\rm dim}\, (\repr_{X})^{\st_{G}(X)}$,
thus
$$d(X) = (- 1)^{n - {\rm dim} \, X}
\sum_{\sur{c \in \chcl [X] \hfill}
{{\rm det} \, (\st_{G} (c)_{\vert X}) \subset \{1\}}}
(- 1)^{|c| - 1},$$
where by ${\rm det} \, (\st_{G} (c)_{\vert X}) \subset \{1\}$
we mean that ${\rm det}\, (\sigma_{|X}) = 1$
for all $\sigma$ in $\st_{G} (c)$.

\begin{cor}The following equalities hold
\begin{equation}
{\rm det}\, (F, H^{n}_{c} (U_{\bar \FF},
q^{\ast} \cL_{\psi})^{G})
=
\prod_{X \in L (\cA) / G}
(\phi ({\rm discr}\, (q_{\vert X})) g (\phi)^{{\rm dim} \,
X})^{d (X)}
\end{equation}
\begin{equation}
{\rm dim}\, (H^{n}_{c} (U_{\bar \FF},
q^{\ast} \cL_{\psi})^{G})_{{ \rm weight}\, k}
=
\sum_{\sur{X \in L(\cA) / G \hfill}
{{\rm dim} \, X = k \hfill}}
d (X).\hfill\qed
\end{equation}
\end{cor}

\begin{prop}The virtual representation $\repr_{X}$
is in fact a representation of
$\st_{G} (X)$, and
thus, in particular,
$d (X) \geq 0$, for any $X$ in $L (\cA)$.
\end{prop}

\begin{proof}Let $\cA_{X}$ be the arrangement in $\AA^{n}$
consisting of all hyperplanes of $\cA$
which contain $X$.
Clearly
$\st_{G} (X)$ acts on $\cA_{X}$.
Let $U_{X}$ be the complement of
$\cA_{X}$ in $\AA^{n}$.
From
Theorem 2.2.5 applied to
$\cA_{X}$ we deduce
that
$\repr_{X}$ is isomorphic, as a virtual representation of
$\st_{G} (X)$, to
$H^{n}_{c} ((U_{X})_{\bar \FF},
q^{\ast} \cL_{\psi})_{{ \rm weight}\, {\rm dim} \, X}$
which is an honest representation.\end{proof}

\begin{cor}Suppose
${\rm dim} \, H^{n}_{c} (U_{\bar \FF},
q^{\ast} \cL_{\psi})^{G} = 1$.
Then there is exactly one
$X_{0}$ in $L (\cA) / G$
such that $d (X_{0}) \not=0$. Moreover, 
$d (X_{0}) = 1$ and
$$
{\rm det}\, (F, H^{n}_{c} (U_{\bar \FF},
q^{\ast} \cL_{\psi})^{G})
=
\phi ({\rm discr}\, (q_{\vert X_{0}})) g (\phi)^{{\rm dim} \,
X_{0}}.
$$
\end{cor}

\begin{proof}The corollary is a direct consequence of Corollary 2.2.6 (2)
and Proposition 2.2.7.\end{proof}

We will see in (4.4) that the hypothesis of Corollary 2.2.8 is satisfied
for a Coxeter arrangement over $\FF$ when $p$ does not divide
$\vert G \vert$.

\bigskip

\subsection{The case when $\chi = \phi$}We keep
the same notations
as in (2.1) and (2.2).

\begin{lem}Suppose that $f$ is the square of a polynomial over
$\FF$ on which $G$ acts by multiplication by
the determinant.
\begin{enumerate}
\item[(1)]
There is  a canonical isomorphism
$$
H^{\cdot}_{c} ((U / G)_{\bar \FF},
f^{\ast}_{G} \cL_{\phi}
\otimes q^{\ast}_{G} \cL_{\psi})
\simeq
H^{\cdot}_{c} (U_{\bar \FF}, q^{\ast} \cL_{\psi})^{\rm det} \, ,
$$
where the superscript
$\rm det$ denotes the part on which $G$
acts by multiplication by
the determinant.
\item[(2)]For $i \not= n$, 
$
H^{i}_{c} ((U / G)_{\bar \FF},
f^{\ast}_{G} \cL_{\phi}
\otimes q^{\ast}_{G} \cL_{\psi}) = 0$.
\item[(3)]The vector space
$H^{n}_{c} ((U / G)_{\bar \FF},
f^{\ast}_{G} \cL_{\phi}
\otimes q^{\ast}_{G} \cL_{\psi})_{{\rm weight} \, n}$
is one dimensional, and the Frobenius
acts on it by multiplication by
$\phi ({\rm discr}\, q)
g (\phi)^{n}$.
\end{enumerate}
\end{lem}

\begin{proof}Let us prove (1).
Let $G^{+} = \{\sigma \in G \bigm| {\rm det} \, \sigma = 1\}$.
We may suppose
$G^{+} \not = G$, otherwise the assertion is trivial.
We factor $\pi : U \rightarrow U /G$
into $\pi_{1} : U \rightarrow U /G^{+}$
and $\pi_{2} : U /G^{+}\rightarrow U /G$.
Note that $\pi_{2}$ is the Kummer cover obtained by
taking the square root
of $f$, so we have canonical isomorphisms
$
f_{G}^{\ast}\cL_{\phi} \simeq
(\pi_{2 \ast} \barql)^{-} \simeq (\pi_{2 \ast} \barql)^{\rm det}.
$ From the
canonical isomorphisms
$(\pi_{2 \ast} \barql)^{\rm det}
\simeq
(\pi_{2 \ast}(\pi_{1 \ast} \barql)^{G^{+}})^{\rm det}
\simeq
((\pi_{2 \ast}\pi_{1 \ast} \barql)^{G^{+}})^{\rm det}
\simeq
(\pi_{\ast} \barql)^{\rm det}$,
we obtain
canonical isomorphisms
\begin{align*}
H^{\cdot}_{c} ((U / G)_{\bar \FF},
f^{\ast}_{G} \cL_{\phi}
\otimes q^{\ast}_{G} \cL_{\psi})
&\simeq
H^{\cdot}_{c} ((U / G)_{\bar \FF},
(\pi_{\ast} \barql)^{\rm det}
\otimes q^{\ast}_{G} \cL_{\psi})\\
&\simeq
H^{\cdot}_{c} ((U / G)_{\bar \FF},
(\pi_{\ast} \barql)
\otimes q^{\ast}_{G} \cL_{\psi})^{\rm det}\\
&\simeq
H^{\cdot}_{c} (U_{\bar \FF}, q^{\ast}_{G} \cL_{\psi})^{\rm det}.
\end{align*}
The assertion (2) is a direct consequence of
(1) and Lemma 2.2.2,
while (3) follows directly
from (1) and Theorem 2.2.5.\end{proof}

\begin{prop}Suppose that the assumptions
of Lemma 2.3.1 are valid and that
$H^{n}_{c} (U_{\bar \FF}, f^{\ast} \cL_{\phi}
\otimes q^{\ast} \cL_{\psi})^{G}$
is of dimension 1. Then
$$
S_{G} (\phi) = (- 1)^{n} \phi ({\rm discr}\, q)
g (\phi)^{n}.
$$
\end{prop}

\begin{proof}It is a direct consequence of the preceding lemma.\end{proof}

\section{Monodromy and critical values for Coxeter
arrangements}
\setcounter{equation}{0}
\subsection{}Let $(V, G, q)$ be a Coxeter arrangement over
$\CC$.
For each reflection hyperplane $H$, we choose a linear form $\ell_{H} : V
\rightarrow \CC$ defining $H$ in $V$ and we set $\delta = (\prod_{H}
\ell_{H})^{2}$. We denote by $N$ the number of reflection hyperplanes
and by $\Delta :
V / G \rightarrow \CC$ the map induced by $\delta$. Thus $\Delta$
is the discriminant of $G$.
We denote by $F_{0}$ the Milnor fiber of $\Delta $ at $0$
and by $Z (T, G)$ the zeta function of local monodromy of $\Delta$ at
$0$, {\it i.e.},
$$
Z (T, G) = \prod_{i} {\rm det} (1 - TM, H^{i} (F_{0}, \CC))^{(-1)^{i + 1}},
$$
where $M$ denotes the monodromy automorphism (see, {\it e.g.},
\cite{Ar-V-G}).

Let $d_{1}, \ldots, d_{n}$ be the degrees of $G$. In \cite{D-L3}
we proved the following result.

\begin{stheorem}For $G$ a finite Coxeter group
we have
$$
\prod_{\cE \, {\rm connected \, subgraph}}
Z (-T, G (\cE))^{(-1)^{\vert \cE \vert}}
=
\prod_{i = 1}^{n} \frac{1 - T^{d_{i}}}{1 -T}
$$
where the product on the left-hand side runs over all connected
subgraphs $\cE$ of the Coxeter diagram of $G$,
$G (\cE)$ denotes the Coxeter group with diagram $\cE$,
and $\vert \cE \vert$ the number of vertices of $\cE$.
\end{stheorem}

The proof of Theorem 3.1 given in \cite{D-L3} depends on a case
by case analysis. A more conceptual proof, but depending on
Macdonald's formula (0.1), is given in \cite{D-L4}.

\medskip

We use now terminology
from
\cite{SGA7} II Exp. XIII - XIV. 
Denote by $\bar \eta_{0}$ the generic geometric point of the
henselization of $\AA^{1}_{\CC}$ at $0$ and by $I_{0, \CC}$
its inertia group ({\it i.e.} the fundamental group of the complement of $0$
in a small disk around $0$).
We denote by $K_{I_{0, \CC}}$ the Grothendieck group
of
finite dimensional vector spaces with
$I_{0, \CC}$-action.
If $\cL$ is an object in $D^{b}_{c} (\GG_{m, \CC}, \CC)$
we denote by
$[\cL_{\bar \eta_{0}}]$ the class of
$\sum (- 1)^{i} [\cH^{i} (\cL)_{\bar \eta_{0}}]$
in $K_{I_{0, \CC}}$
and we set
$[\cL_{\bar \eta_{\infty}}] =
[{\rm inv}^{\ast}(\cL)_{\bar \eta_{0}}]$, where $\rm inv$
is the morphism $x \mapsto x^{-1}$.
If a finite group $G$ acts on $\cL$ we denote by
$[\cL^{G}_{\bar \eta_{0}}]$ the class of
$\sum (- 1)^{i} [\cH^{i} (\cL)^{G}_{\bar \eta_{0}}]$
and we define similarly $[\cL^{G}_{\bar \eta_{\infty}}]$.
For any character $\chi : I_{0, \CC} \rightarrow \CC^{\times}$ we denote
by $V_{\chi}$ the class in
$K_{I_{0, \CC}}$ of the rank one object with action given by $\chi$
and for any natural number $m \geq 1$,
we set
$V_{m} = [(\pi_{m \ast} \CC)_{\bar \eta_{0}}]$,
where $\pi_{m} : \GG_{m} \rightarrow \GG_{m}$ is given by $x \mapsto x^{m}$.
Of course, we have $V_{m} = \sum_{\chi^{m} = 1} V_{\chi}$.
We set
$\bar M_{G} = (-1)^{n - 1} [(R\psi_{\Delta} (\CC))_{0}]$.
Here $(R\psi_{\Delta} (\CC))_{0}$ is the stalk  at zero of the complex
of nearby cycles with respect to $\Delta$.

One can rephrase Theorem 3.1. as follows.
\newtheorem*{TTT}{Theorem 3.1$'$}
\begin{TTT}For $G$ a finite Coxeter group,
the following equality in $K_{I_{0, \CC}}$ holds,
$$
\sum_{\cE \, {\rm connected \, subgraph}} \bar M_{G (\cE)}
=
V_{\phi}
\otimes (\sum_{i= 1}^{n} (V_{d_{i}}
- V_{1})),
$$
where $\phi$ is the unique character of order 2.
\end{TTT}

\subsection{}Let $U$ be the complement in $V$ of $\delta = 0$,
let $B$ (resp. $B_{0}$) be the intersection of $U$ with the locus
of $q= 1$ (resp. $q = 0$).
The following proposition is proved in \cite{D-L4} by a detailed analysis of
the geometry of a nice compactification of $U$.
Here and in the following we denote by 
$q^{N}_{\vert U}$ the restriction of the map defined by $q$
on $U$ raised to $N$-th power.

\begin{prop}We have the following equalities in $K_{I_{0, \CC}}$,
\begin{enumerate}
\item[(1)]
$[(R \delta_{\vert B_{0} !} \CC)^{G}_{\bar \eta_{\infty}}]
=
[(R \delta_{\vert B !} \CC)^{G}_{\bar \eta_{\infty}}]$.
\item[(2)]$[(R \delta_{\vert B !} \CC)^{G}_{\bar \eta_{\infty}}]
+
[(Rq^{N}_{\vert U !} \CC)^{G}_{\bar \eta_{0}}]
=
[(R\psi_{\Delta} (\CC))_{0}]$.
\item[(3)]$[(R \delta_{\vert B !} \CC)^{G}_{\bar \eta_{0}}]
=
(- 1)^n \sum_{\sur{\cE \, {\rm connected \, subgraph}}{{G (\cE) \not= G}}} \bar M_{G (\cE)}$.
\item[(4)]There exist
$\bar a, \bar b$ in $\ZZ$ with $\bar a + \bar b = (-1)^{n - 1}$ such that
$[(Rq^{N}_{\vert U !} \CC)^{G}_{\bar \eta_{0}}]
=
(\bar a - \bar b) V_{N} + \bar b V_{2N}$
and such that
$[(Rq_{\vert U !} \delta^{\ast} \cL_{\chi})^{G}_{\bar \eta_{0}}]
=
\bar a V_{\chi^{N}} + \bar b V_{\phi \chi^{N}}$,
for any character $\chi$ of $I_{0, \CC}$, where $\cL_{\chi}$ denotes
the Kummer sheaf associated to $\chi$.
\end{enumerate}
\end{prop}

\begin{cor}The equality
$$
[(R \delta_{\vert B !} \CC)^{G}_{\bar \eta_{0}}]
-
[(R \delta_{\vert B !} \CC)^{G}_{\bar \eta_{\infty}}]
-
[(Rq^{N}_{\vert U !} \CC)^{G}_{\bar \eta_{0}}]
=
(-1)^{n}(V_{\phi}
\otimes \sum_{i= 1}^{n} (V_{d_{i}}
- V_{1}))
$$
holds
in
$K_{I_{0, \CC}}$.
\end{cor}

\begin{proof}The corollary follows directly from (2), (3) and Theorem 3.1$'$.\end{proof}

\medskip

Actually, the values of $\bar a$ and $\bar b$
are determined by the following proposition, which is proved in \cite{D-L4},
but not needed in the present paper.

\begin{prop}The following formula holds,
$$[(Rq_{\vert U !} \CC)^{G}_{\bar \eta_{0}}]
= (-1)^{n - 1} V_{\phi^{n + N}}
= \bar a V_1 + \bar b V_{\phi}.$$
\end{prop}

\subsection{}If $f : X \rightarrow \GG_{m, \CC}$ is a morphism of
schemes of finite type over $\CC$,
we set
$a (f) = \prod_{s \in \GG_{m, \CC}} s^{a_{s} (Rf_{!} \CC)}$
where $a_{s}$ is the drop of the rank of $Rf_{!} (\CC)$ at $s$
({\it i.e.} the alternating sum of the drops of the rank
of its cohomology sheaves, {\it cf.} \cite{La} 3.1.5.2.
We denote by $\Delta_{B}$
the function
$B / G \rightarrow \GG_{m, \CC}$
induced by $\delta$.
We will use the following result which is proven in \cite{D-L4}
as an easy consequence of Macdonald's formula (0.1),
by which one calculates the critical value (0.5) of
$\Delta_B$.

\setcounter{stheorem}{2}
\begin{stheorem}The following formula holds,
$$
a (\Delta_{B})
=
\Bigl(\kappa
\frac{
\prod_{i}d_{i}^{d_{i}}}
{N^{N}}\Bigr)^{(-1)^{n}},
$$
where $\kappa$ is as in (0.2).
\end{stheorem}

\section{Proof of the Main Theorem}
\setcounter{equation}{0}
\subsection{Determinants and monodromy}Let $U$
be a scheme of finite type
over the finite field $\FF$. We assume that
the characteristic $p$ of $\FF$ is not equal to 2 and we
fix a prime $\ell$ distinct from $p$.
We denote by $\bar \FF$ an algebraic closure of $\FF$
and we set
$U_{\bar \FF} = U \otimes_{\FF} \bar \FF$.
For any object $\cF$ in $D_{c}^{b} (U, \barql)$
(the derived category of ``bounded complexes"
of $\barql$-sheaves with constructible cohomology) we denote
$$
\varepsilon_{0} (U / \FF, \cF) := {\rm det} (-F,
H_c^{\cdot}
(U_{\bar \FF}, \cF))^{-1}.
$$
Let $f : U \rightarrow \GG_{m, \FF}$ be a morphism.
We set
$a (f) = \prod_{s \in \GG_{m, \bar \FF}} s^{a_{s} (Rf_{!} \barql)} \in \FF$,
where $a_{s}$ denotes the total drop at $s$ of a complex of
sheaves
({\it cf.} \cite{La} 3.1.5.2). We also denote by $\bar \eta_{0}$ the generic
geometric point of the henselization of $\AA^{1}_{\bar \FF}$ at $0$
and by $I_{0}$ (or $I_{0, \bar \FF}$) the inertia group.
We denote by $K_{I_{0}}$ the Grothendieck group
of finite dimensional $\barql$-vector spaces with continuous action
of $I_{0}$ defined on a finite extension of $\QQ_{\ell}$.
For $\cL$ an object in $D_{c}^{b} (\GG_{m, \FF}, \barql)$
we define as in (3.1) objects
$[\cL_{\bar \eta_{0}}]$ and
$[\cL_{\bar \eta_{\infty}}]$ in $K_{I_{0}}$. Also if a finite
group $G$ acts on $\cL$ we define similarly $[\cL^{G}_{\bar \eta_{0}}]$ and
$[\cL^{G}_{\bar \eta_{\infty}}]$.
For $\chi : I_{0} \rightarrow \barql^{\times}$ a continuous character,
we denote by
$V_{\chi}$ the associated object in $K_{I_{0}}$. For $N \geq 1 $ an
integer, we define $V_{N}$ in $K_{I_{0}}$ as in (3.1).
We will denote by $\phi$ the character of order 2 of
$I_{0}$.
For any character $\chi : \FF^{\times} \rightarrow \barql^{\times}$
we denote by
$\cL_{\chi}$ the corresponding Kummer sheaf on $\GG_{m, \FF}$.
We can also view $\chi$ as a continuous
character $I_{0} \rightarrow \barql^{\times}$ which we still denote by
$\chi$.

Suppose $[(Rf_{!} \barql)_{\bar \eta_{0}}] = \sum_{N \in \NN^{\times}}
\alpha_{N}
V_{N}$
and
$[(Rf_{!} \barql)_{\bar \eta_{\infty}}] = \sum_{N \in \NN^{\times}}
\beta_{N}
V_{N}$
with $\alpha_{N}$ and $\beta_{N}$ in $\ZZ$. (Here $\NN^{\times}$ denotes
$\NN \setminus \{0\}$.)
We define the rational number $b (f)$ as
$$
b (f) := \Bigl(\prod_{N \in \NN^{\times}} N^{N \alpha_{N}}\Bigr)
\Bigl(\prod_{N \in \NN^{\times}}  (-N)^{- N \beta_{N}}\Bigr).
$$
For any character $\chi : \FF^{\times} \rightarrow \barql^{\times}$ we set
$$
G (f, \chi) = \prod_{N \in \NN^{\times}}
\Bigl(g (\chi^{N})^{\alpha_{N}}
g (\chi^{-N})^{\beta_{N}}\Bigr)
$$
and
$$
\tilde G (f, \chi) = \prod_{N \in \NN^{\times}}
\Bigl(g (\chi^{N})^{\alpha_{N} - \beta_{N}}\Bigr)
$$
where $g (\chi)$ is the
Gauss sum
$g (\chi) := -\sum_{x \in \FF^{\times}} \chi (x) \psi (x)$
for $\psi$ a fixed non trivial additive character.

\begin{theorem}Assume that $\alpha_{N} = \beta_{N} = 0$
when $p$ divides $N$ and that the sheaves
$R^{i}f_{!} \barql$ are tame at $0$ and $\infty$
for every $i$. Then, for every multiplicative character
$\chi$ of $\FF^{\times}$,
\begin{align*}
\varepsilon_{0} (U / \FF, f^{\ast} \cL_{\chi})
&=
\varepsilon_{0} (U / \FF, \barql) \,
\chi (\frac{a (f)}{b(f)}) \, G (f, \chi)\\
&=\varepsilon_{0} (U / \FF, \barql)
\, \vert \FF \vert^{\sum_{\chi^{N} \not= 1} \beta_{N}}
\, \chi (\frac{a (f)}{|b(f)|}) \, \tilde G (f, \chi).
\end{align*}
\end{theorem}

\begin{proof}The theorem follows directly from \cite{D-L2} Proposition
2.4.1, which is based on Laumon's product formula \cite{La} 3.2.1.1, using
\cite{La} 3.1.5.4 (iv),
and
$g (\chi^{- 1})
= \vert \FF \vert \chi (-1) g (\chi)^{-1}$ if $\chi \not= 1$.
For more details, see \cite{D-G} Proposition 3.2.1. A much stronger result
is contained in \cite{Sa}.\end{proof}

\noindent {\it Remark}. --- The advantage
of the second equality upon the first is that
$|b (f)|$ and
$\tilde G (f, \chi)$ only depend on
$[(Rf_{!} \barql)_{\bar \eta_{0}}] -
[(Rf_{!} \barql)_{\bar \eta_{\infty}}]$.

\bigskip

Let $\cF$ be
a constructible $\ell$-adic sheaf on a scheme $U$
(of finite type) over $\FF$
and suppose 
a finite group
$G$ acts on $U$ and $\cF$.
We set
$$
\varepsilon_{0}^{G} (U / \FF, \cF) := {\rm det} (-F, H^{\cdot}_c
(U_{\bar \FF}, \cF)^{G})^{-1}.
$$
\medskip

\noindent {\it Notation}. --- For $\chi$ any multiplicative character of
$\FF^{\times}$ and for $F$ in
$K_{I_{0, \CC}}$ of the form
$F = \sum_{N \in \NN^{\times}}
\gamma_{N} V_{N}$, with $\gamma_{N}$ in $\ZZ$, we set
$$
\tilde G (F, \chi) = \prod_{N} g (\chi^{N})^{\gamma_{N}}
$$
and
$$
b (F) = \prod_{N} N^{N \gamma_{N}}.
$$
It is assumed here that a non trivial additive character
has been fixed. Note that $\tilde G$ and $b$ are multiplicative in $F$.

\medskip

\begin{sloppypar}
\begin{lem}
\begin{enumerate}
\item[(1)]If $N \in \NN^{\times}$ is even, 
$V_{N} \otimes V_{\phi} = V_{N},$
$\tilde G (V_{N} \otimes V_{\phi}, \chi)
=
\tilde G (V_{N}, \chi)$ and
$b (V_{N} \otimes V_{\phi}) = b (V_{N}).$
\item[(2)]If $N \in \NN^{\times}$ is odd, 
$V_{N} \otimes V_{\phi} = V_{2N} - V_{N}$,
$\tilde G (V_{N} \otimes V_{\phi}, \chi)
= \frac{\chi (4^{N})}{g (\phi)} \,
\tilde G (V_{N}, \phi \chi)$
and $b (V_{N} \otimes V_{\phi}) = 4^{N}
b (V_{N}).$ 
\end{enumerate}
\end{lem}
\end{sloppypar}

\begin{proof}Everything is clear except perhaps the second relation in (2) which
follows directly from the
Hasse-Davenport
formula ({\it cf.} \cite{E1}) and the first one.\end{proof}

\subsection{}Fix an integer $d \geq 1$.
Let $R \subset \CC$ be a Dedekind domain
containing the  $d$-th roots of unity, with fraction field $K$
of finite degree over $\QQ$.
We denote by 
$\bar K = \bar \QQ$ the algebraic closure of $K$ in $\CC$.
Let $A = (V, G, q)$ be a Coxeter arrangement over $R$.
For $\got P$  a maximal ideal of $R$ we denote by
$({\overline V}_{\got P}, {\overline G}_{\got P}, {\overline q}_{\got
P})$
the corresponding data over $k_{\got P} := R / \got P$.
When no confusion can arise, we will omit the index $\got P$.
Choose linear forms
$\ell_{H}$ over $R$ which define the hyperplanes of $\cA_G$.
Let $\got P$ be a maximal ideal of $R$ which is relatively prime
to
$d \, \vert G \vert$, such that the linear forms
$\ell_{H}$ are not zero modulo $\got P$ and such that
${\overline q}_{\got P}$ is a nondegenerate bilinear form.
(In particular $\got P$ is relatively prime  to 2 and to each degree
$d_{i}$ of $G$, because $\vert G \vert$ is even and equal to the product
of the $d_{i}$'s.)
Thus the conclusions of Propositions 1.2.2, 1.3.1, and 1.3.3 hold.
In this situation there is no ambiguity for the notations,
$\delta_{\got P}$, $U_{\got P}$,
$B_{\got P}$, $B_{0 \got P}$, etc, {\it cf.} the notation in 3.1 and 3.2,
which refer to objects
attached to the reduction of $A$ mod $\got P$, which are also the
reduction mod $\got P$ of the corresponding objects over $R$.
We denote by $I_{0, \bar K}$ the inertia group at 0
of $\GG_{m, \bar K}$ and we fix a character
$\chi : I_{0, \bar K} \rightarrow \barql^{\times}$ of order
$d$. The Kummer sheaf $\cL_{\chi}$ on
$\GG_{m, \bar K}$ is obtained by base change from a sheaf, which we still
denote
by $\cL_{\chi}$, on $\GG_{m, R}$. 
The character $\chi$ induces a character
$\chi_{\got P} : I_{0, \bar k_{\got P}} \rightarrow \barql^{\times}$
(here $\bar k_{\got P}$ denotes an algebraic closure of
$k_{\got P}$) which is associated to a character,
still denoted by the same symbol, $\chi_{\got P} : k_{\got P}^{\times}
\rightarrow \barql^{\times}$.
For $s$ an element of $K$ we write
$\chi_{\got P} (s)$ for the value of $\chi_{\got P}$
on the residue class of $s$
if $s$ is a unit in  $R_{\got P}$ and
we set
$\chi_{\got P} (s) = 1$
otherwise.

\begin{prop}Let $\got P$ be a maximal ideal of $R$
satisfying the above conditions. Then
${a (\Delta_{B})} / {b (F)}$ and $\kappa$ are units in $R_{\got P}$
and 
\begin{equation}
\begin{split}
&\chi_{\got P}
\Bigl(\frac{a (\Delta_{B})}{b (F)}\Bigr) \,
\tilde G (F, \chi_{\got P})
\, =\\
&\qquad \qquad C
\chi_{\got P} (\kappa)^{(-1)^{n}}
\prod_{1 \leq i \leq n}
\Bigl[\frac{g ((\phi \chi_{\got P})^{d_{i}})}{g (\phi \chi_{\got
P})}\Bigr]^{(-1)^{n}}
\frac{g (\chi_{\got P}^{N})^{\bar a} g (\phi \chi_{\got P}^{N})^{\bar
b}}{g (\phi)^{\bar b}}
\end{split}
\end{equation}
for
$F = [(R \delta_{\vert B !} \CC)^{G}_{\bar \eta_{0}}]
-
[(R \delta_{\vert B !} \CC)^{G}_{\bar \eta_{\infty}}]$
in $K_{I_{0, \CC}}$
(it follows from   \S\kern .15em 3
and Lemma 4.1.2 that $F$
has the required form),
with $C = \Bigl[g (\phi)^{(\sum_{d_{i} \, {\rm odd}} 1) - n}\Bigr]^{(-1)^{n -
1}}$.
Here, as usual, $\phi$ denotes the character of order 2.
\end{prop}

\begin{proof}It follows from Proposition 1.3.3 that $\kappa$
is a unit in $R_{\got P}$. By Corollary
3.2.2,
$$
F - E = (-1)^{n}(V_{\phi}
\otimes \sum_{i= 1}^{n} (V_{d_{i}}
- V_{1}))
$$
with $E = [(Rq^{N}_{\vert U !} \CC)^{G}_{\bar \eta_{0}}]$,
and by 3.2.1 (4),
$E = \bar a  V_{N} + \bar b (V_{2 N} - V_{N})$.
So by Lemma 4.1.2, we obtain
that $b (F - E)$ is a unit in $R_{\got P}$
and
$$
\frac{\tilde G (F - E, \chi_{\got P})}{\chi_{\got P} (b (F - E))}
=
C \prod_{1 \leq i \leq n}
\Bigl[\frac{g ((\phi \chi_{\got P})^{d_{i}})}{g (\phi \chi_{\got
P})}\Bigr]^{(-1)^{n}}
\chi_{\got P} \Bigl(\prod_{1 \leq i \leq n} d_{i}^{d_{i}}\Bigr)^{(-1)^{n
- 1}},
$$
while $b (E) N^{(-1)^{n}N}$
is a unit in $R_{\got P}$ and
$$
\frac{\tilde G (E, \chi_{\got P})}{\chi_{\got P} \Bigl(b (E)
N^{(-1)^{n}N}\Bigr)}
=
\frac{
g (\chi_{\got P}^{N})^{\bar a}
g (\phi \chi_{\got P}^{N})^{\bar b}}
{g (\phi)^{\bar b}}.
$$
(To verify this last equation when $N$ is even, one has again to use
the Hasse-Davenport formula as in the proof of 4.1.2.)
Now the  formula  follows directly from the multiplicativity of
$\tilde G$ and $b$ and from
Theorem 3.3.\end{proof}

\setcounter{equation}{0}
\begin{prop}For almost all $\got P$ in ${\rm Spec \,}
R$ ({\it i.e.} for all but a finite number),
\begin{equation}
\varepsilon_{0}^{G} (B_{\got P} / k_{\got P}, \delta^{\ast}_{\got P}
\cL_{\chi_{\got P}})
=
\varepsilon_{0}^{G} (B_{\got P} / k_{\got P}, \barql)
\vert k_{\got P} \vert^{\sum_{\chi^{j} \not= 1} \beta_{j}}
\chi_{\got P}
(\frac{a (\Delta_{B})}{b (F)})
\tilde G (F, \chi_{\got P})
\end{equation}
and
\begin{equation}
\varepsilon_{0}^{G} (B_{0 \got P} / k_{\got P}, \delta^{\ast}_{\got P}
\cL_{\chi_{\got P}})
=
\varepsilon_{0}^{G} (B_{0 \got P} / k_{\got P}, \barql)
\vert k_{\got P} \vert^{\sum_{\chi^{j} \not= 1} \beta_{j}},
\end{equation}
where the $\beta_j$ are given by
$[(R \delta_{\vert B !}\CC)^G_{\bar \eta_{\infty}}] = \sum_j \beta_j V_j$.
\end{prop}

\begin{proof}The first assertion is a direct consequence of Theorem 4.1.1,
using standard constructibility, comparison and base change theorems.
For the second assertion, note that because $\delta$ and
$q$ are homogeneous
the sheaves
$R^{i} \delta_{\vert B_{0} !} \CC$ are lisse on $\GG_{m}$. Thus
$\vert b (\Delta_{\vert B_{0} / G}) \vert = 1$
(with $b$ defined as in 4.1)
and
$[(R \delta_{\vert B_{0} !} \CC)^{G}_{\bar \eta_{0}}] =
[(R \delta_{\vert B_{0} !} \CC)^{G}_{\bar \eta_{\infty}}]$.
Now the  result follows similarly from
Theorem 4.1.1 and
Proposition 3.2.1 (1) (that
$[(R \delta_{\vert B_{0} !} \CC)^{G}_{\bar \eta_{\infty}}]
=
[(R \delta_{\vert B !} \CC)^{G}_{\bar \eta_{\infty}}]$).\end{proof}

\medskip

\subsection{Good reduction}Let $A = (V, G, q)$ be a
liftable
Coxeter arrangement over $\FF$. We choose
defining linear forms
for each reflection hyperplane.
Let $d \geq 1$ be an integer and
$\chi : \FF^{\times} \rightarrow \barql^{\times}$ a character
of order $d$. There exists  a
discrete valuation ring $T \subset \CC$,
with residue field $\FF$
and fraction field of finite degree
over $\QQ$, containing the $d$-th roots of unity,
and a lifting $A_{T}$ of $A$ over $T$.
We denote by $K$ the fraction field of $T$, by $\bar K$ its
algebraic closure in $\CC$
and by $\cL_{\chi}$ the Kummer sheaf on $\GG_{m, T}$ associated to
$\chi$. We assume that $p$ does not divide
$\vert  G \vert$ and we choose defining linear forms
with coefficients in $T$
for each reflection hyperplane, which reduce
to those already chosen in $\FF$. 
We denote by $\delta_{T}$, $U_{T}$, $B_{T}$, etc, the
data associated to $A_{T}$.

\begin{lem}If $p$ does not divide
$\vert  G \vert$,
the specialization morphism
$$
H^{i}_{c} (B  \otimes \bar \FF, \delta^{\ast}
\cL_{\chi})
\longrightarrow
H^{i}_{c} (B_{T} \otimes \bar K, \delta_{T}^{\ast}
\cL_{\chi})
$$
is an isomorphism for all $i$.
\end{lem}

\begin{proof}The arrangement associated to $A_{T}$ has good reduction by
Proposition 1.3.1. Thus
the $T$-scheme $U_{T}$ admits a canonical smooth compactification $\bar U_{T}$
over $T$ such that $\bar U_{T} \setminus U_{T}$ has normal crossings over $T$,
and such that $\bar U_T$ is
equipped with a morphism $\pi : \bar U_{T} \rightarrow
\PP^{n}_{T}$, where $\PP^{n}_{T}$ is the projective closure of $V_{T}$,
extending the inclusion of $U_{T}$ in $V_{T}$. The compactification
$\bar U_{T}$ and the morphism
$\pi$ are obtained by blowing up successively the union of the
strict transforms of the strata of dimension $i$
of the projective arrangement
associated to $A_{T}$ for $i$ increasing from $0$ to $n - 2$
({\it cf.} \cite{L} \S \kern .15em 7). We now define
$\bar B_{T}$ to be the strict transform
in $\bar U_{T}$ of the closure of $B_{T}$ in $\PP^{n}_{T}$.
It follows from Proposition 1.3.3 that
$\bar B_{T}$ intersects  the strata of
$\bar U_{T} \setminus U_{T}$ transversally. Hence $\bar B_{T}$ is smooth and proper over $T$ and
$\bar B_{T} \setminus  B_{T}$ has normal crossings over $T$. The assertion now
follows from
\cite{SGA41/2} Th. finitude Appendice 1.3.3 and 2.4.\end{proof}

\begin{cor}
Let
$A = (V, G, q)$ be a
liftable
Coxeter arrangement over $\FF$. If
$p$ does not divide
$\vert  G \vert$, then
$$
\sum (- 1)^i {\rm dim} \, H^i_c (B \otimes \bar \FF,
\delta^{\ast} \cL_{\chi})^G = (- 1)^{n - 1}.
$$
\end{cor}

\begin{proof}The isomorphism in Lemma 4.3.1 being $G$-equivariant, we only have
to prove the analogous result for Coxeter arrangements over $\CC$.
In this case
\begin{align*}
\sum (- 1)^i {\rm dim} \, H^i_c (B,
\delta^{\ast} \cL_{\chi})^G &=
\sum (- 1)^i {\rm dim} \, H^i_c (B,
\CC)^G\\
&=\vert G \vert^{- 1} \chi (B, \CC).
\end{align*}
But in \cite{D-L4} we proved that
$\chi (B, \CC) = (- 1)^{n - 1} \vert
G \vert$ for complex Coxeter arrangements.
(The proof is by an easy induction in a more general setting,
showing that $(- 1)^{n - 1} \chi (B, \CC)$
equals the number of Weyl chambers.)\end{proof}

\begin{lem}
If $p$ does not divide
$\vert  G \vert$, the sheaves
$R^{i}q_{T \vert U_T !} \delta_{T}^{\ast} \cL_{\chi}$ are lisse on
$\GG_{m, T}$ for all $i$.
\end{lem}

\begin{proof}The polynomials associated to
$\delta_{T}$ and $q_{T}$ being homogeneous it is enough to
prove that the specialization
morphisms associated to the specialization of the generic point of
$t = 1$ to the special point are isomorphims, which is precisely the
content of Lemma 4.3.1.\end{proof}

\begin{cor}Let
$A = (V, G, q)$ be a
liftable
Coxeter arrangement over $\FF$. If
$p$ does not divide
$\vert  G \vert$, the sheaves
$R^{i}q_{\vert U !} \delta^{\ast} \cL_{\chi}$ are lisse and tame on
$\GG_{m, \FF}$ for all $i$.
\end{cor}

\begin{proof}This is a direct consequence of Lemma 4.3.3, see, {\it e.g.},
\cite{D-L1} (4.1).\end{proof}

\medskip

Localizing the ring of algebraic integers in $K$ by inverting a finite number of elements,
we obtain a Dedekind ring $R \subset T$ and a Coxeter arrangement
$A_R$ over $R$, such that $A_T$ is obtained from
$A_R$ by extension of scalars.
Let ${\got P}_0 \in {\rm Spec} \, R$ be the intersection of the maximal ideal
of $T$ with $R$. Then ${\got P} = {\got P}_0$ satisfies the conditions stated at the beginning of
4.2.

\begin{prop}For ${\got P} = {\got P}_0$ or for
${\got P}$ any other maximal ideal of $R$ satisfying the conditions
stated
at the beginning of 4.2, the following holds:
\begin{enumerate}
\item[(1)]The representation of
${\rm Gal} (\bar K \, \vert \, K)$
on
${\rm det} (H^{\cdot}_c (B_T \otimes \bar K, \delta^{\ast}
\cL_{\chi})^G)^{-1}$
is unramified at $\got P$
and the action of the geometric Frobenius relative to $\got P$ is
given by
$- \varepsilon_{0}^{G} (B_{\got P} / k_{\got P}, \delta^{\ast}_{\got P}
\cL_{\chi_{\got P}})$.
\item[(2)]The representation of
${\rm Gal} (\bar K \, \vert \, K)$
on
${\rm det} (H^{\cdot}_c (B_{0 T} \otimes \bar K, \delta^{\ast}
\cL_{\chi})^G)^{-1}$
is unramified at $\got P$
and the action of the geometric Frobenius relative to $\got P$ is
given by
$\varepsilon_{0}^{G} (B_{0 \got P} / k_{\got P}, \delta^{\ast}_{\got P}
\cL_{\chi_{\got P}})$.
Moreover the Euler characteristic of 
$R \Gamma_{c} (B_{0 \got P}  \otimes  \bar k_{\got P},
\delta^{\ast}_{\got P}
\cL_{\chi_{\got P}})$ is zero.
\end{enumerate}
\end{prop}

\begin{proof}The first statement in (1) follows directly from Lemma
4.3.1 and smooth base change. The second statement in (1) follows now
from Corollary 4.3.2.
For (2), observe that by
Lemma 4.3.3 and  the Leray spectral sequence the analogous statement holds
for the Galois module
${\rm det} (H^{\cdot}_c (U_T \setminus  B_{0 T} \otimes \bar K, \delta^{\ast}
\cL_{\chi})^G)^{-1}$, by the same argument as at the end
of the proof of Lemma 4.3.1, using
the compactification $\PP^1_T$ of $\GG_{m, T}$.
So it is enough to prove
the analogue statement for
${\rm det} (H^{\cdot}_c (U_T  \otimes \bar K, \delta^{\ast}
\cL_{\chi})^G)^{-1}$, which is clear if one considers the compactification
used in the proof of Lemma 4.3.1. The assertion about
the Euler characteristic follows in the same way.\end{proof}

\renewcommand{\theequation}{\thesubsection.\arabic{equation}}
\setcounter{equation}{0}
\subsection{Proof of the Main Theorem}Let $A = (V, G, q)$ be a
Coxeter arrangement over $\FF$. We assume that $p$ does not divide
$\vert G \vert$.
Let $d \geq 1$ be an integer and
$\chi : \FF^{\times} \rightarrow \barql^{\times}$ a character
of order $d$. By Proposition 1.6, the Coxeter arrangement $A$ is liftable.
We use the notation and material of 4.3. It follows from
Proposition 1.3.3 that
$\kappa \in T \setminus {\got P}_0$. It is well known (see
\cite{B} Ch.V, \S\kern .15em 6 Th\'{e}or\`{e}me 1) that
$\sum_{i = 1}^{n} (d_i - 1) = N$. This implies 
that the right-hand side of
4.2.1 (1) may be expressed in terms of a Jacobi sum Hecke character
({\it cf.} \cite{SGA41/2} Sommes trig.
\S\kern .15em 6).
By Proposition
4.3.5, Proposition 4.2.2 and Cebotarev density,
we deduce that relations (1) and (2)
in Proposition 4.2.2 are also valid for ${\got P}_0$. Together with
Proposition 4.2.1 this gives
\begin{equation}
\begin{split}
& \varepsilon_{0}^{G} (B / \FF, \delta^{\ast}
\cL_{\chi})
= \varepsilon_{0}^{G} (B / \FF, \barql)
\\
& \qquad \qquad
\cdot \vert \FF \vert^{\sum_{\chi^{j} \not= 1} \beta_{j}}
C \chi (\kappa)^{(-1)^{n}} \prod_{1 \leq i \leq n}
\Bigl[\frac{g ((\phi \chi)^{d_{i}})}{g (\phi \chi)}\Bigr]^{(-1)^{n}}
\frac{g (\chi^{N})^{\bar a} g (\phi \chi^{N})^{\bar b}}{g (\phi)^{\bar b}}
\end{split}
\end{equation}
and
\begin{equation}
\varepsilon_{0}^{G} (B_0 / \FF, \delta^{\ast}
\cL_{\chi})
=
\varepsilon_{0}^{G} (B_0 / \FF, \barql)
\vert \FF \vert^{\sum_{\chi^{j} \not= 1} \beta_{j}}.
\end{equation}
Here we have implicitly used
Proposition 1.5 to identify the degrees of
$A$ and $A_T$.

By Lemma 4.3.3 the sheaves
$R^{i}q_{T \vert U_T !} \delta_{T}^{\ast} \cL_{\chi}$ are lisse on
$\GG_{m, T}$, hence we deduce from Proposition 3.2.1 (4)
that
$[(Rq_{\vert U !} \delta^{\ast} \cL_{\chi})^G_{\bar \eta_0}]
=
\bar a V_{\chi^N} + \bar b V_{\phi \chi^N}$.
The sheaves
$(R^{i}q_{\vert U !} \delta^{\ast} \cL_{\chi})^G$ being lisse and tame on
$\GG_{m, \FF}$ by Corollary 4.3.4, this implies, by the Leray spectral sequence,
the relation
\begin{equation}
\varepsilon_{0}^{G} (U \setminus B_0 / \, \FF, \delta^{\ast}
\cL_{\chi} \otimes q^{\ast} \cL_{\psi})
=
g (\chi^N)^{\bar a} g (\phi \chi^N)^{\bar b}
\varepsilon_{0}^{G} (B / \FF, \delta^{\ast}
\cL_{\chi})^{-1}.
\end{equation}
(To verify this one uses the structure of tame irreducible lisse sheaves on
$\GG_{m, \FF}$ (see, {\it e.g.}, \cite{La} p.198), or, alternatively,
one can use \cite{D-G} 8.1.4 and 3.1.4,(4) and \cite{La} 3.1.5.5.)
Since we have
\begin{equation}
\begin{split}
& \varepsilon_{0}^{G} (U/ \FF, \delta^{\ast}
\cL_{\chi} \otimes q^{\ast} \cL_{\psi}) = \\
& \qquad \qquad \qquad
\varepsilon_{0}^{G} (B_0 / \FF, \delta^{\ast}
\cL_{\chi}) \, 
\varepsilon_{0}^{G} (U \setminus B_0 / \, \FF, \delta^{\ast}
\cL_{\chi} \otimes q^{\ast} \cL_{\psi}),
\end{split}
\end{equation}
from equations (4.4.1) - (4.4.4), we obtain the relation
\begin{equation}
\begin{split}
& \varepsilon_{0}^{G} (U/ \FF, \delta^{\ast}
\cL_{\chi} \otimes q^{\ast} \cL_{\psi})
=\frac{\varepsilon_{0}^{G} (B_0/ \FF, \barql)}{\varepsilon_{0}^{G} (B/
\FF, \barql)} \\
& \qquad \qquad
\cdot C^{- 1}
\chi (\kappa)^{(-1)^{n - 1}} g (\phi)^{\bar b}
\prod_{1 \leq i \leq n}
\Bigl[\frac{g ((\phi \chi)^{d_{i}})}{g (\phi \chi)}\Bigr]^{(-1)^{n - 1}}
\, .
\end{split}
\end{equation}
So,
$$
{\rm det} \, (- F, H^{\cdot}_c
(U_{\bar \FF}, \delta^{\ast}
\cL_{\chi} \otimes q^{\ast} \cL_{\psi})^G)
=
A \Bigl[\chi (\kappa)
\prod_{1 \leq i \leq n}
\frac{g ((\phi \chi)^{d_{i}})}{g (\phi \chi)}\Bigr]^{(-1)^{n}},
$$
with $A$ independent of the character $\chi$.

Now note that by Proposition 2.1.2 and Corollary 4.3.4,
$H^{i}_c (U_{\bar \FF}, \delta^{\ast} \cL_{\chi} \otimes q^{\ast}
\cL_{\psi})^G = 0$
for $i \not= n$. On the other hand,
$$
\chi \Bigl(R \Gamma_c (U_{\bar \FF}, \delta^{\ast} \cL_{\chi} \otimes q^{\ast}
\cL_{\psi})^G \Bigr)
=
\chi \Bigl(R \Gamma_c (\GG_{m, \bar \FF}, (Rq_{\vert U !}
\delta^{\ast} \cL_{\chi})^G \otimes \cL_{\psi})\Bigr),
$$
because of Proposition 4.3.5 (2).
Since the $(R^iq_{\vert U !}
\delta^{\ast} \cL_{\chi})^G$ are lisse and tame on $\GG_{m, \FF}$
we obtain
$$
\chi \Bigl(R \Gamma_c (\GG_{m, \bar \FF}, (Rq_{\vert U !}
\delta^{\ast} \cL_{\chi})^G \otimes \cL_{\psi})\Bigr)
=
- \sum_i (-1)^i {\rm dim} \, H^i_c (B \otimes \bar \FF, \delta^{\ast}
\cL_{\chi})^G,
$$
see, {\it e.g.}, \cite{K} 4.8.2.
By Corollary 4.3.2, we deduce that
\begin{equation}
{\rm dim} \, H^{n}_c (U_{\bar \FF}, \delta^{\ast} \cL_{\chi} \otimes q^{\ast}
\cL_{\psi})^G = 1.
\end{equation}
Hence we obtain by the Grothendieck trace formula, 
\begin{equation}
\begin{split}
S_G (\chi) &= {\rm Tr} (F,
H^{\cdot}_c (U_{\bar \FF}, \delta^{\ast} \cL_{\chi} \otimes
q^{\ast} \cL_{\psi})^G),\\
&=
A' \chi (\kappa)
\prod_{1 \leq i \leq n}
\frac{g ((\phi \chi)^{d_{i}})}{g (\phi \chi)},
\end{split}
\end{equation}
with $A'$ independent of the character $\chi$.
To find the value of $A'$ we take $\chi = \phi$,
so $A' = S_G (\phi) \phi (\kappa)$
and by Proposition 2.3.2 we obtain
\begin{equation}
A' = \phi (\kappa) (- 1)^n \phi ({\rm discr} q) g (\phi)^{n}.
\end{equation}
The theorem now follows from (4.4.7) and (4.4.8).\hfill$\qed$

\subsection{Remark}A slightly different way of organizing the proof of
the Main Theorem is by using material in \cite{A2} or \cite{L-S} (or the results of
\cite{Sa-T}). This material provides an analogue of Theorem 4.1.1
for the determinant of period integrals. In this way Macdonald's formula (0.4),
which yields the period $\int_{B (\RR)} \Delta^{s + 1 / 2} \frac{dx}{dq}$,
directly implies an analogous expression for
$\varepsilon_0^G (B / \FF, \delta^{\ast} \cL_{\chi}) \, 
\varepsilon_0^G (B / \FF, \barql)^{- 1}$, when $\chi$ is generic and
$p \gg 0$, because the last expression and the period are given by similar
formulas (compare Theorem 4.1.1 with {\it loc. cit.}).
To obtain a full proof of the Main Theorem
by this approach, one has to do similar calculations as in (4.4),
but now one needs Proposition 3.2.1 (1) and Proposition 3.2.3.

\subsection{Remark}That the analogy between our Main Theorem
and Macdonald's formula (0.4) is no coincidence, can be explained partially
by the following conjecture of Deligne (\cite{De} \S\kern .15em 8.9):
the period of a rank one motive over a number field is completely determined
by the Frobenius action. Using this conjecture one directly obtains
from Macdonald's formula an expression for $\varepsilon_0^G (B / \FF, \delta^{\ast} \cL_{\chi})$
when $\chi$ is generic and
$p \gg 0$. Note that to apply the conjecture one first has to get rid of the additive character.
(For this reduction one needs again Proposition 3.2.1 (1) and Proposition 3.2.3.)
\bibliographystyle{plain}
\bibliography{trans}

\end{document}